\documentclass[a4paper]{article}

\usepackage{amsmath, amsthm}
\usepackage{amssymb} 
\usepackage[latin1]{inputenc} 
\usepackage[english]{babel}
\usepackage{graphicx}
\usepackage{tikz}
\usepackage{xspace} 
\usepackage{tabls}  
\usepackage{float} 
\usepackage{enumerate}
\usepackage{lscape}

\usepackage{fancyhdr}
\title{\ttl}
\fancyheadoffset{1cm}
\fancyfootoffset{.4cm}
\newcommand{\ben}{\begin{enumerate}}
\newcommand{\een}{\end{enumerate}}
\newcommand{\bit}{\begin{itemize}}
\newcommand{\eit}{\end{itemize}}
\newcommand{\bce}{\begin{center}}   
\newcommand{\ece}{\end{center}}
\newcommand{\ig}{\includegraphics}
\newcommand{\btab}{\begin{tabular}}  
\newcommand{\etab}{\end{tabular}}
\newcommand{\bma}{\begin{matrix}}
\newcommand{\ema}{\end{matrix}}
\newcommand{\barr}{\begin{array}}
\newcommand{\earr}{\end{array}}

\newcommand{\bfig}[3]{
\bce
\begin{figure}[H]
 #1
\vspace{.2cm}
\caption{#2}
\label{#3}
\end{figure}
\ece
}

\newcommand{\matr}[1]{\left( \bma #1 \ema \right)}

\newcommand{\bal}[1]{\begin{align*}
#1
\end{align*}}

\newcommand{\mip}[2]{
\begin{minipage}[h!]{#1 \textwidth}
  #2 
\end{minipage}
}

\newcommand{\bb}{\mathbb} 

\newcommand{\R}{\ensuremath{\bb{R}}\xspace}

\newcommand{\N}{\ensuremath{\bb{N}}\xspace}

\usepackage[hidelinks]{hyperref}

\newcommand*
	\ttl{Bürgi's ''Kunstweg'' -- geometric approach}

\newcommand*\klss{}                                       

\newcommand*\lhrr{C.\ Riedweg}       

\newcommand*\jahr{}       

\newcommand*\planung{}       

\author{Christian Riedweg}

\begin{document}

\maketitle
\begin{abstract}
	Placing regular $4n$-sided polygons correctly,  equations between sine and sums of sines  show up -- exactly these equations are used in Bürgi's method to approximate sines.
%
%
\end{abstract}



%

$ $

$ $

Menso Folkerts, Dieter Launert and Andreas Thom
describe in
 \cite{flt_sinus}  a method for calculating sines. This method was discovered by Jost Bürgi  (1552-1632) who called it ''Kunstweg''.
Folkerts, Launert and Thom  prove that this method really works i.e. they show  the convergence of this method.

One of the steps in this proof 
is  to show that 
for $\alpha = 90^\circ / n$ (where $n\in \N$)
\[
	\matr{
		\sin(\alpha)
		\\ 
		\sin(2\alpha)
		\\ 
		\sin(3\alpha)
		\\ 
		\sin(4\alpha)
		\\ 
		\vdots 
		\\ 
		\sin\big((n-1)\alpha\big) 
		\\
		 1
	}
	\quad
	\text{is an eigenvector of}
	\quad
	\matr{ 1 & 1 & 1 & 1 & \ldots & 1 & 1/2
			\\
			1 & 2 & 2 & 2 &\ldots& 2 & 2/2
			\\
			1 & 2 & 3& 3& \ldots & 3 & 3/2
			\\
			1 & 2 & 3& 4& \ldots & 4 & 4/2
			\\
			\vdots &\vdots&\vdots&\vdots& \ddots &\vdots&\vdots
			\\
			1 & 2 & 3 & 4 & \ldots  & n-1 & (n-1)/2 
			\\
			1 & 2 & 3 & 4 & \ldots & n-1 & n/2 
	} 
	\in \R^{n\times n}.
\]
I.e. one has to show that for some $\lambda$ the following equations hold:
\bal{
	\lambda \sin(\alpha) 
		&= 
			\sin(\alpha)
			+ 
			\sin(2\alpha)
			+ 
			\ldots 
			+ 
			\sin\big((n-1)\alpha\big) 
			+
			1/2
	\\
	 & \ldots 
	\\
	 \lambda \sin(j\alpha)&= 
		\sin(\alpha) 
		+
		2\sin(2\alpha)
		+ 
		\ldots 
		+ 
		j\sin(j\alpha) 
		+ 
		j\sin\big((j+1)\alpha\big) 
		+ 
		\ldots 
		+
		j\sin\big((n-1)\alpha\big) 
		+ j/2
	\\
	 & \ldots 
	\\
	\lambda \quad 
		&=
		\sin(\alpha) 
		+
		2\sin(2\alpha)
		+
		3\sin(3\alpha)
		+ 
		\ldots 
		+ 
		(n-1)\sin\big((n-1)\alpha\big) 
		+  n/2.
}

\vspace{.3cm}
Subtracting subsequently the preceeding equation, this system of equations is equivalent to

\mip{.9}
	{
	\bal
		{
		\lambda \sin(\alpha) 
			&= 
				\sin(\alpha)
				+ 
				\sin(2\alpha)
				+ 
				\ldots 
				+ 
				\sin\big((n-1)\alpha\big) 
				+
				 1/2
		\\
		\lambda \sin(2\alpha)
			-
			\lambda \sin(\alpha) 
			&=
			\sin(2\alpha)
			+ 
			\sin(3\alpha)
			+ 
			\ldots 
			+ 
			\sin\big((n-1)\alpha\big) 
			+ 1/2
		\\
		& 
		 \ldots 
		\\
		\lambda \sin\big( j\alpha\big)
			-
			\lambda \sin\big( (j-1)\alpha\big)
			&=
			\sin(j\alpha)
			+ 
			\sin\big((j+1)\alpha\big)
			+ 
			\ldots 
			+ 
			\sin\big((n-1)\alpha\big) 
			+ 1/2
		\\
		& 
		 \ldots 
		\\
		\lambda 
		 - 
			\lambda \sin ((n-1)\alpha)
		 &= 
			1/2.	 
	}
}
\mip{.05}{
	\vspace{.5cm}
	$(\star)$
}

\vspace{.5cm}

In  section~\ref{sec:geo_app}, we  develop this system of equations $(\star)$ purely geometrically (hence we give a geometric proof that the above vector is an eigenvector).
\\
In section~\ref{sec:details}, we discuss 
in more details how these equations are derived from regular $4n$-sided polygons.
\\
For  more information on Bürgi's method and further references, see \cite{flt_sinus}.

\newpage

\section{Geometric approach to Bürgi's ''Kunstweg''}
	\label{sec:geo_app}

Let $n\in\N$.
We place regular convex $4n$-sided polygons along to each other: the second  not exactly on top of the first but ''one step to the left''. The third is placed on top of the second but ''two steps to the left'', and so on, as visualized below for $12$-sided polygons:

\bfig
	{
	\ig
		[scale=1.6]
		{figure1}
	}
	{
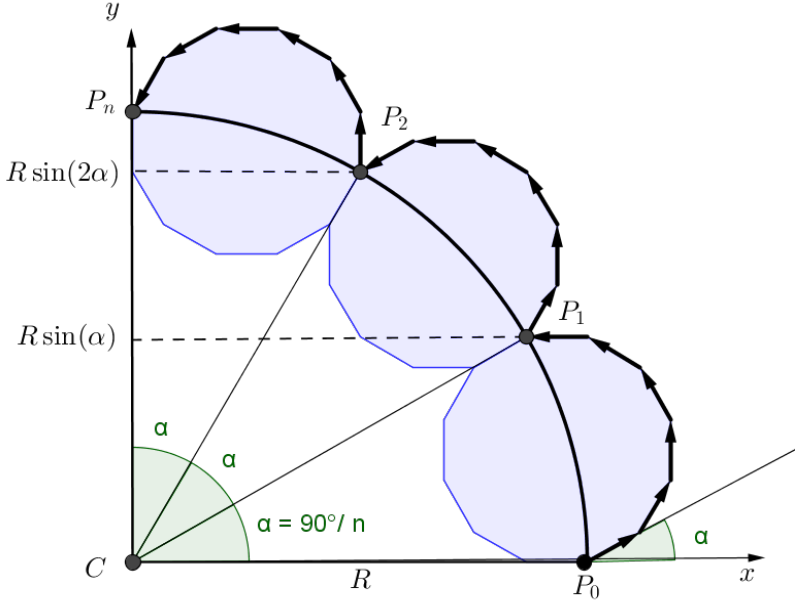
	For the $4n$-sided polygon is $\alpha = 90^\circ /n$ (more details in the next section). Hence we need $n$ such polygons to fill an angle of $90^\circ$ at $C$.
	}
	{fig1}

\vspace{.2cm}
The $y$-coordinates of the points $P_1$, $P_2$, \ldots,  $P_n$ are
\[
	R\sin(\alpha), 
	\quad
	R\sin(2\alpha),
	\quad
	\ldots,
	\quad
	R\sin((n-1)\alpha),
	\quad
	R. 
\]
Thus, the $y$-coordinate variation from point $P_{j-1}$ to $P_j$ is
\[
	R\sin\big( j\alpha\big)
			-
			R \sin\big( (j-1)\alpha\big).
\]

On the other hand, 
following the indicated vectors (all of lenght $l$) from $P_{j-1}$ to $P_j$    we get 
a sum of sines (more details -- see next section). 
\\
Together:
\bal{
	R\sin(\alpha)
		&=
		 2l\sin(\alpha) + 2l\sin(2\alpha) + \ldots + 2l\sin((n-1)\alpha) + l,
	\\
	&\vdots
	\\
		R\sin\big( j\alpha\big)
			-
			R \sin\big( (j-1)\alpha\big)
			&=
			2l\sin(j\alpha)
			+ 
			2l\sin\big((j+1)\alpha\big)
			+ 
			\ldots 
			+ 
			2l	\sin\big((n-1)\alpha\big) 
			+ l
	\\
	&\vdots
	\\
	R - 
		R\sin((n-1)\alpha) 
			&=
		 l.
}
These are exactly the equations $(\star)$ that show up in Thom's proof in \cite{flt_sinus} (after dividing by $2l$ we get $\lambda = R/2l$).

\newpage

\section{Regular convex 4n-sided polygons}
\label{sec:details}
Let $n\in \N$; we consider regular convex $4n$-sided polygons (visualized are $12$-sided polygons). Let all these polygons have a side parallel to the $x$-axis.

\bfig
	{\hspace{.5cm}
	\ig
		[scale=1.2]
		{figure2}
	}
	{
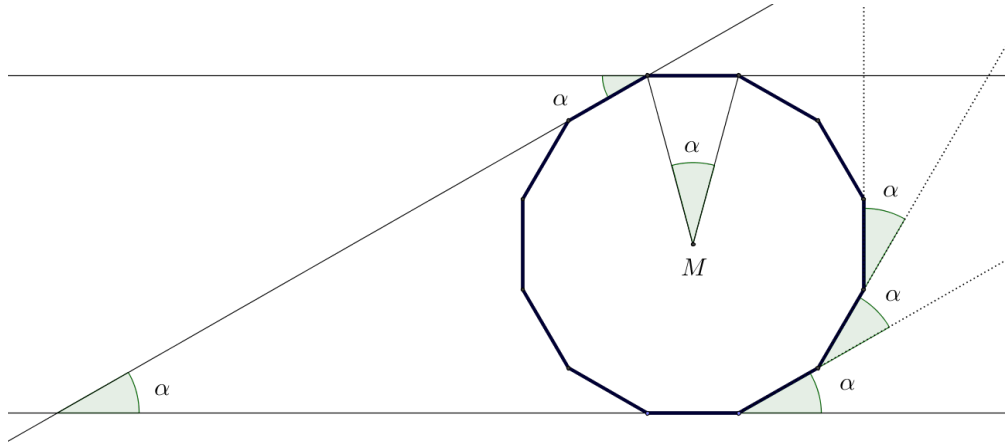
	$M$ is the center of the regular polygon. The angle $\alpha$ occurs at various places.
	For  $4n$-sided polygons we get
	$\alpha = \frac{360^\circ}{4n} = \frac{90^\circ}n.$
	}
	{fig2}

Now we consider different paths  along such polygons. These $n$ paths all feature a lenght  of \textbf{half  of the polygon's circumference}, i.e. each path consits of $2n$ sides of the polygon (for regular $4n$-sided polygons)
\\
As we are interested  in sines only,  the $y$-coordinate of the  vectors along the paths has to be considered. Thus, horizontal vectors do not contribute at all and some pairs of vectors  cancel each others $y$-coordinate -- both cases are indicated by a dashed line below.

\vspace{.5cm}
\bit
	\item {The \textbf{1. polygon}, first path}
	\bfig
		{\hspace{2cm}
			\ig
				[scale=1.1]
				{figure3}
		}
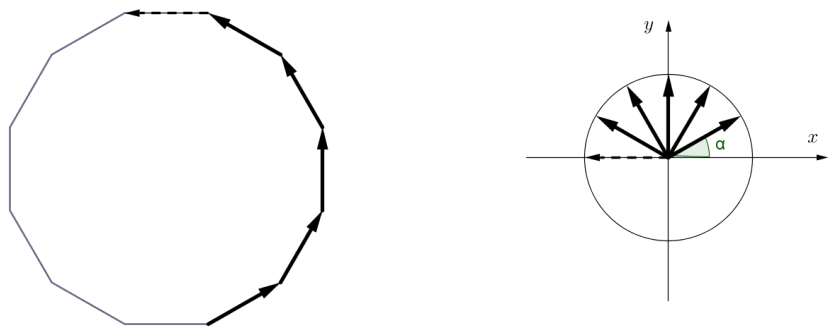
		{The $2n$-vectors in the left picture are visualized in the coordinate plane (right picture). All these vectors are of the same lenght, say $l$. In Figure~\ref{fig2} we see that the angle between the $x$-axis and the first vector is $\alpha$, it is  $2\alpha$ for the second vector, $3\alpha$ for the third, and so on. 
		}
		{fig3}
	
	Following these vectors, the change of the $y$-coordinate equals  
	\[
		l\cdot\sin(\alpha) + l\cdot\sin(2\alpha) + \ldots + l\cdot \sin(2n\alpha).
	\]
		
		With $\alpha = \frac{90^\circ}n$, the considered angles are
		\bce
			 $\frac{90^\circ}n$, 
			 \quad
			 $2\cdot \frac{90^\circ}n$, 
			 \quad
			 \ldots, 
			 \quad
			 $n\alpha=90^\circ$, 
			 \quad
			 \ldots, 
			 \quad
			 $2n\alpha = 180^\circ$.
		\ece
		Since
		 $\sin(90^\circ)=1$,  $\sin(180^\circ)=0$ and  $\sin(\varphi)=\sin(180^\circ -\varphi)$ we can reduce this to
	
		\[
			2l\sin\left(\alpha\right) 
			+
			2l\sin(2\alpha)
			+
			\ldots
			+
			2l\sin\big((n-1)\cdot \alpha\big)
			+
			l.
		\]
			
	\vspace{1cm}
	
	\item {The
			\textbf{$j.$ polygon}, where $1<j\le n$  (visualized is $j=2$).}
	\bfig
		{\hspace{2cm}
			\ig
			[scale=.6]
			{figure5}
			\hspace{2cm}
			\ig
			[scale=.6]
			{figure4}
		}
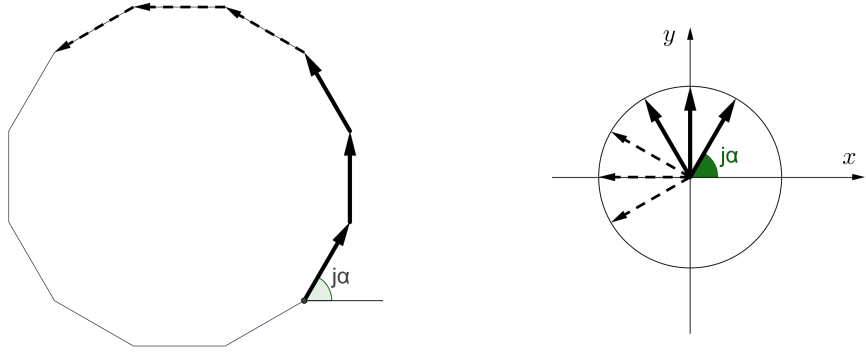
		{Let the angle of the first vector  be $j\alpha$. Again, the length of the vectors is $l$. 
		}
		{fig4}
	Here, the angles between the $x$-axis and the considered vectors  are
	\[
		j\alpha, \quad
		(j+1)\alpha, \quad
		\ldots, \quad
		(j+2n-1)\alpha.
	\]

		\vspace{.2cm}
		Hence the $y$-coordinate changes by
		\[
			l\cdot\sin(j\alpha) + l\cdot\sin\big((j+1)\alpha\big)+\ldots + l\cdot\sin((j+2n-1)\alpha).
		\]
		We can simplify this sum:
		\bal{
			\text{for $j<n$:} \qquad & 2l\sin(j\alpha)
			+
			2l\sin\big((j+1)\alpha\big)
			+
			\ldots
			+
			2l\sin\big((n-1)\alpha\big)
			+
			l
			\\
			\text{for $j=n$:} \qquad & l.
		}

	\vspace{.5cm}
	\newpage
	\item
	\textbf{All polygons together:}
	For symmetrical reasons,  the rays starting in the points $P_j$ to the lower left (as pictured below) intersect all in one point -- say $C$. Let $C$ be the origin.
	\bfig{\hspace{2cm}
			\ig
				[scale=1.6]
				{figure1}
		}
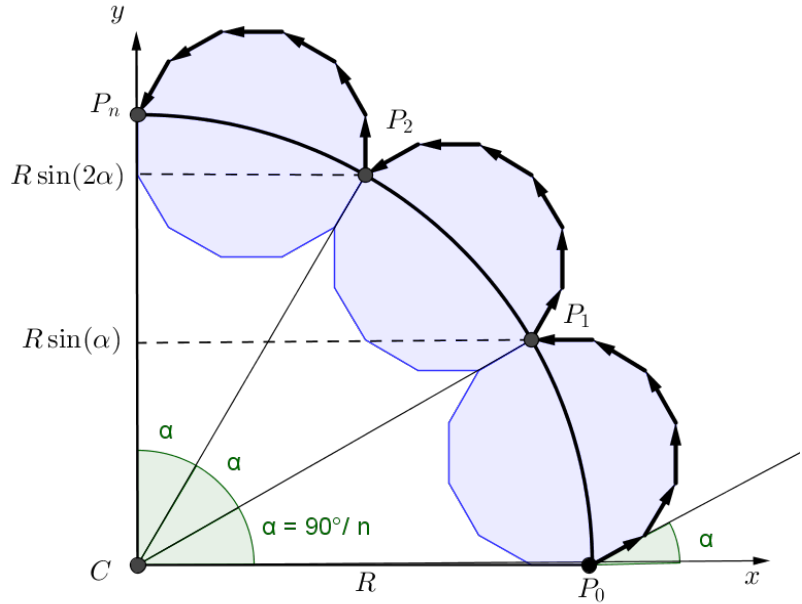
		{The points $P_0$, $P_1$, \ldots, $P_n$ are on a circle (with center $C$) of some radius, say $R$. 	Note that  $P_n$ is on the $y$-axis since $\alpha = 90^\circ / n$. }
		{}

	Starting in $P_0$ and following the indicated path we pass the points $P_1$, $P_2$, \ldots, $P_n$. By symmetry, these points are on a common circle with center $C$. Let $R$ be the radius of this circle.
	Hence the $y$-coordinate of $P_j$ is  $R\sin(j\alpha)$. On the other hand, this $y$-coordinate is equal to the $y$-coordinate of $P_{j-1}$ plus the $y$-coordinate change from $P_{j-1}$ to $P_j$ as seen above.  Together we get the system of equations $(\star)$ of page 1 (after dividing by $l$ and with $\lambda=R/l$):
	for $P_1$:
	\[
		R \sin(\alpha) = 2l\sin(\alpha) + 2l \sin\big(2\alpha\big) + \ldots + 2l\sin\big((n-1)\alpha\big) + l,
	\]
	for $P_j$, where $1<j<n$:
	\[
		R\sin(j\alpha)
		-
		R\sin\big( (j-1)\alpha\big)
		=
		2l\sin(j\alpha) + 2l \sin\big((j+1)\alpha\big) + \ldots + 2l\sin\big((n-1)\alpha\big) + l
	\]
	and for $P_n$:
	\[
		R-R\sin\big((n-1)\alpha\big) = l.
	\]

\eit

\end{document}